\documentclass[onefignum,onetabnum]{siamart251216}

\newcommand{\norm}[1]{\left\lVert#1\right\rVert}

\newcommand{\A}[0]{\mathcal A}

\newcommand{\M}[0]{\mathcal M}

\newcommand{\etaweak}[0]{\eta_{\text{weak}}}
\newcommand{\etastrong}[0]{\eta_{\text{strong}}}

\newcommand{\SoftwarePackage}[1]{\textsc{#1}} %
\newcommand{\belos}{\SoftwarePackage{Belos}}
\newcommand{\muelu}{\SoftwarePackage{MueLu}}

\newcommand{\teko}{\SoftwarePackage{Teko}}
\newcommand{\trilinos}{\SoftwarePackage{Trilinos}}

\newcommand{\klu}{\SoftwarePackage{KLU2}}
\newcommand{\amesos}{\SoftwarePackage{Amesos2}}
\newcommand{\ifpack}{\SoftwarePackage{Ifpack2}}
\newcommand{\zoltan}{\SoftwarePackage{Zoltan2}}
\newcommand{\frosch}{\SoftwarePackage{FROSch}}
\newcommand{\aria}{\SoftwarePackage{Aria}}

\let\mychp\section
\let\mysec\subsection
\let\mysubsec\subsubsection

\usepackage{booktabs}
\usepackage{lipsum}
\usepackage{amsfonts}
\usepackage{graphicx}
\usepackage{epstopdf}
\usepackage{algorithmic}
\usepackage{cleveref}
\usepackage{tikz}
\usepackage{amsmath}
\usepackage{siunitx}
\usepackage{tikz}
\usepackage{tikz-3dplot}
\usepackage{pgfplots}
\usepackage{pgfplotstable}
\usetikzlibrary{shapes,arrows,calc,fit,backgrounds,shapes.multipart,decorations.pathreplacing}
\usepgfplotslibrary{groupplots}
\pgfplotsset{compat = newest}
\usetikzlibrary{intersections}
\usepgfplotslibrary{fillbetween}

\graphicspath{{figures/}}
\ifpdf
  \DeclareGraphicsExtensions{.eps,.pdf,.png,.jpg}
\else
  \DeclareGraphicsExtensions{.eps}
\fi

\newsiamremark{remark}{Remark}
\newsiamremark{hypothesis}{Hypothesis}
\crefname{hypothesis}{Hypothesis}{Hypotheses}
\newsiamthm{claim}{Claim}
\newsiamremark{fact}{Fact}
\crefname{fact}{Fact}{Facts}

\headers{Preconditioning for Thermal Batteries}{M. Phillips}

\title{Multi-physics Preconditioning for Thermally Activated Batteries}

\author{Malachi Phillips\thanks{Sandia National Laboratories, Albuquerque, NM (\email{malphil@sandia.gov}).}}

\usepackage{amsopn}

\ifpdf
\hypersetup{
  pdftitle={Preconditioning for Thermal Batteries},
  pdfauthor={M. Phillips}
}
\fi

\begin{document}

\maketitle

\begin{abstract}
  Thermal batteries, also known as molten-salt batteries, are single-use reserve power systems activated by pyrotechnic heat generation, which transitions the solid electrolyte into a molten state.
The simulation of these batteries relies on multiphysics modeling to evaluate performance and behavior under various conditions.
This paper presents advancements in scalable preconditioning strategies for the Thermally Activated Battery Simulator (TABS) tool, enabling efficient solutions to the coupled electrochemical systems that dominate computational costs in thermal battery simulations.
We propose a hierarchical block Gauss-Seidel preconditioner implemented through the Teko package in Trilinos, which effectively addresses the challenges posed by tightly coupled physics, including charge transport, porous flow, and species diffusion.
The preconditioner leverages scalable subblock solvers, including smoothed aggregation algebraic multigrid (SA-AMG) methods and domain-decomposition techniques, to achieve robust convergence and parallel scalability.
Strong and weak scaling studies demonstrate the solver's ability to handle problem sizes up to 51.3 million degrees of freedom on 2048 processors, achieving near sub-second setup and solve times for the end-to-end electrochemical solve.
These advancements significantly improve the computational efficiency and turnaround time of thermal battery simulations, paving the way for higher-resolution models and enabling the transition from 2D axisymmetric to full 3D simulations.

\end{abstract}

\mychp{Introduction}
\label{chp:intro}
Thermal batteries, also referred to as molten-salt batteries, are single-use primary reserve power systems that utilize a solid electrolyte at ambient temperatures, which transitions to a molten state upon activation \cite{Crompton1982, Guidotti2006}.
Activation is achieved by igniting pyrotechnic pellets, which generate heat to liquefy the electrolyte, enabling the battery to function until the electrolyte re-solidifies or the reactants are exhausted.
Due to their single-use design, multiphysics modeling is essential for developing and evaluating the performance of these batteries.

Recent advancements in TABS, or the Thermally Activated Battery Simulator Tool, \cite{voskuilen2021multi,roberts2023development,roberts2018tabs,voskuilen2023coupled} have focused on improving the fidelity and usability of multiphysics models for thermal batteries.
Voskuilen et al. \cite{voskuilen2021multi} explored high-fidelity full models for thermal battery simulations, which are currently constrained to 2D axisymmetric representations due to limitations to linear solvers for the coupled electrochemical system.
For example, the reported calculation wall time for the 2D axisymmetric full model consisted of nearly a full day of compute, with the majority of the time spent in the linear solver for the electrochemical system.
The ultimate objective of this research is to enable full 3D simulations of these detailed multiphysics battery models, which would provide a more comprehensive understanding of battery behavior and performance, especially under the effects of mechanical deformation.
Achieving this goal requires the implementation of efficient linear solvers, which are critically dependent on effective preconditioning strategies to ensure computational feasibility and scalability.

We consider here the linear solution to the coupled electrochemical system as described in \Cref{sec:electrochemical}, a tightly-coupled multiphysics system requiring monolithic preconditioning approaches.
Matrix snapshots are generated from simulations using the SIERRA Multimechanics Module: Aria \cite{ariaV526}, henceforth referred to as \aria{} for brevity.
Preconditioning strategies to the linear system rely on the \trilinos{} software project \cite{Trilinos,Heroux2005a,mayr2025trilinos},
which focuses on a package-level design for addressing the various computational requirements for computational science and engineering applications.

Prior to this work, the state of multiphysics preconditioning in \aria{} included only sparse direct solvers and the domain-decomposition techniques later discussed in \Cref{sec:domain-decomp}.
As the latter methods prove ineffective for the complex multiphysics system described in \Cref{sec:electrochemical}, sparse direct solvers were the only viable path to solving the coupled linear systems.
As a result, both strong and weak scalability were poor.
This work seeks to address these concerns through the design of block-based, physics-aware preconditioners through the \teko{} package in \trilinos{} \cite{cyr2016teko}.

The use of block-based preconditioning techniques for multiphysics problems has been considered in many different contexts, including Navier--Stokes \cite{elman2002block,elman2003parallel,cyr2012stabilization}, 
magnetohydrodynamics \cite{ohm2024scalable,ma2016robust,phillips2014block}, and ice sheet dynamics \cite{isaac2015solution,brown2012composable,helanow2025theoretical}, to name a few.
More relevant to our immediate application, Fang et al. \cite{fang2019parallel} consider the use of block Gauss-Seidel approaches to preconditioning multiphysics systems arising in lithium-ion cells.
In the current work, we consider block-based preconditioners for the linear systems arising from the coupled electrochemical system encountered in the simulation of thermal batteries in the TABS project.
We employ a hierarchical block Gauss-Seidel approach to effectively decouple interactions between the charge transport equations and other terms, including continuity and multi-component Stefan-Maxwell diffusion, in the preconditioner.

The structure of this work is as follows.
The coupled physics in the thermal battery electrochemical system are described in \Cref{sec:electrochemical}.
\Cref{chp:preconditioning} details two approaches to multiphysics preconditioning, including one-level Schwarz methods in \Cref{sec:domain-decomp}
and physics-based block preconditioning in \Cref{sec:block-preconditioning}.
Results are shown in \Cref{chp:results}, which is further split into evaluating the subblock level solvers in \Cref{sec:subblock-solvers},
followed by the end-to-end electrochemical solve in \Cref{sec:end-to-end-solver}.
After this, conclusions follow in \Cref{chp:conclusion}.

\mysec{Coupled Electrochemical Physics in Thermal Batteries}
\label{sec:electrochemical}

\begin{figure}
\centering
\begin{tikzpicture}

\def\width{4}
\def\height{0.5}
\def\collectorHeight{0.2}
\def\insulation{0.5}
\def\can{0.5}
\def\ambient{0.5}
\def\numLayers{20}

\draw[fill=red!50] (0,0) rectangle (\width,\height) node[midway] {Heat Pellet};

\draw[] (0,\height) rectangle (\width,\height+\collectorHeight);
\draw[fill=blue!50] (0,\height+\collectorHeight) rectangle (\width,2*\height+\collectorHeight) node[midway] {Anode};
\draw[fill=gray!50] (0,2*\height+\collectorHeight) rectangle (\width,3*\height+\collectorHeight) node[midway] {Separator};
\draw[fill=yellow!50] (0,3*\height+\collectorHeight) rectangle (\width,4*\height+\collectorHeight) node[midway] {Cathode};
\draw[] (0,4*\height+\collectorHeight) rectangle (\width,4*\height+2*\collectorHeight);

\draw[decorate,decoration={brace,amplitude=10pt,mirror,raise=4pt},yshift=0pt]
(0.7*\width+0.1, \height+0.1) -- (0.7*\width+0.1, 4*\height+2*\collectorHeight-0.1) node[midway,right, xshift=0.35cm] {\( N \)};

\draw[fill=red!50] (0,4*\height+2*\collectorHeight) rectangle (\width,4*\height+2*\collectorHeight+\height) node[midway] {Heat Pellet};

\draw[fill=gray!30] (\width,0) rectangle (\width+\insulation,4*\height+2*\collectorHeight+\height) node[midway,rotate=90] {Insulation};
\draw[fill=black!20] (\width+\insulation,0) rectangle (\width+\insulation+\can,4*\height+2*\collectorHeight+\height) node[midway,rotate=90] {Can};

\draw[] (\width+\insulation+\can,0) rectangle (\width+\insulation+\can+\ambient,4*\height+2*\collectorHeight+\height) node[midway,rotate=90] {Ambient};

\draw[dashed] (0,0) -- (0,4*\height+2*\collectorHeight+\height);
\node[rotate=90] at (-0.2,2*\height+\collectorHeight+\height/2) {Axis};

\draw[<-] (\width/2, \height/2+\height-\collectorHeight) -- (\width+1.5, \height/2+\height-\collectorHeight) node[right] {Collector};
\draw[<-] (\width/2, 4*\height+\collectorHeight+\height/2-0.5*\collectorHeight) -- (\width+1.5, 4*\height+\collectorHeight+\height/2-0.5*\collectorHeight) node[right] {Collector};

\end{tikzpicture}
\caption{2D axisymmetric simulation domain for multi-physics simulations (not to scale).
Note that the collector, anode, separator, and cathode layers are repeated $N=20$ times.}
\label{fig:thermal_battery_stack}
\end{figure}
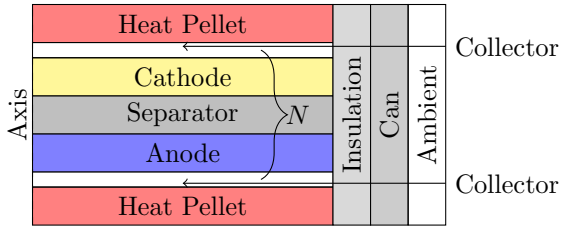

The electrochemical system is a tightly coupled multiphysics problem that demands carefully designed physics-based monolithic preconditioning strategies.
As such, we present a brief overview of the various physics that contribute to the system and provide a few insights into the construction of the block preconditioner later in \Cref{sec:block-preconditioning}.
The following description is a summarization of the electrochemical system for the thermal batteries as described by Voskuilen et al. \cite{voskuilen2021multi}.

We start our discussion of the electrochemical system with evaluating the transport equations associated with the solid and liquid-phase voltages.
The solid-phase voltage, $\Phi_s$, is described via Ohm's law
\begin{equation}\label{eq:solid-phase-voltage}
    \nabla\cdot\left(-\sigma \nabla\Phi_s\right)=S_e,
\end{equation}
where $S_e$ is a source term based on the production of electrons based on the electrochemistry described by \cite{voskuilen2021multi} and derived in \cite{Guidotti2006,masset2008thermal,guidotti2008thermally,tomczuk1982phase,wen1981chemical}.
The electrical conductivity is represented by $\sigma$, which varies by as much as 10 orders-of-magnitude across the simulation domain shown in \Cref{fig:thermal_battery_stack}.
For example, $\sigma \approx 10^{-4}$ \unit{\siemens\per\meter} in the separator, while the conductivity is as large as $\sigma \approx 10^{6}$ \unit{\siemens\per\meter} in the anode.
The collector, anode, separator, and cathode layers of the battery are repeated multiple times throughout the simulation domain, as shown in \Cref{fig:thermal_battery_stack}.
As such, the solid-phase voltage is a highly variable coefficient Poisson equation, which is amenable to the algebraic multigrid (AMG) techniques described in \Cref{sec:sa-amg}.
A homogeneous Dirichlet boundary condition is applied to one of the collectors, with the remaining boundaries observing inhomogeneous Neumann boundary conditions.

The liquid-phase voltage, $\Phi_l$, is described by a more complex relationship,
\begin{equation}\label{eq:liquid-phase-voltage}
    \underbrace{\tau_I\dfrac{\partial\Phi_l}{\partial t}}_{\text{Advective}} + \underbrace{\nabla\cdot\left(F\sum_j \overrightarrow{D}_j z_j\right)}_{\text{Diffusive}}=-S_e,
\end{equation}
where $\overrightarrow{D}_j$ represents the ion diffusive flux, $z_j$ represents the ion valence for each liquid-phase species $j$, and $F$ is Faraday's constant.
The inductive term on the left-hand side of the liquid voltage equation is added for numerical stability, with $\tau_I$ being a small constant.
As a result, the diffusive term in \Cref{eq:liquid-phase-voltage} dominates.
This causes the liquid-phase voltages to approximately obey a diffusion-dominated advection-diffusion system.
Remarkably, the same AMG techniques that are used for the solid-phase voltage may be applied here.
No explicit boundary conditions are applied, as the liquid-phase voltages are coupled to the solid-phase voltage
through the Butler-Volmer coupling
\begin{equation}\label{eq:butler-volmer}
    i = i_0 \left[\exp{\left(\dfrac{\beta Fz\eta}{RT}\right)} - \exp{\left(\dfrac{-(1-\beta)Fz\eta}{RT}\right)}\right],
\end{equation}
where $i$ is the current, $i_0$ is the current density, and $\eta$ is the overpotential.
The form presented in \Cref{eq:butler-volmer} is itself a re-arrangement of the net reaction rate at each electrode surface across each reaction $i$ \cite{bessler2007influence},
\begin{equation}\label{eq:reaction-rate}
  \omega_i = A_i \exp{\left(\dfrac{-E_{a,i}}{RT}\right)}\exp{\left(\dfrac{-\beta F\Delta(\Phi z)_i}{RT}\right)}\left(C_f-C_r\exp{\left(\dfrac{\Delta\mu_i}{RT}\right)}\right),
\end{equation}
where $\beta=0.5$, $\Delta\mu_i$ is the difference in chemical potential, $C_f$ and $C_r$ are the forward and reverse activity coefficients, respectively, and $\Delta(\Phi z)_i$ is the difference in electrical potential energy between the reactants and products.

With the charge transport equations established in the two preceding paragraphs,
we now focus on the porous flow and species transport equations associated with the electrochemical system.
The porosity, or non-solid volume fraction, is treated as a constant
\begin{equation}\label{eq:porosity}
    \phi = \left(1 - \dfrac{V_{solid}}{V}\right) \dfrac{V}{V^{\prime}},
\end{equation}
where $V$ is the volume, $V_{solid}$ the solid-phase volume, and $V^{\prime}$ the volume after any mechanical deformation.
Given the lack of any mechanical deformation, $V^{\prime}=V$ in this work.
The intrinsic permeability, $\kappa_i$, is defined via the Carman-Kozeny model
\begin{equation}\label{eq:carman-kozeny}
    \kappa_i = \dfrac{1}{2 S_v^2 \tau^2} \dfrac{\phi^3}{(1-\phi)^2},
\end{equation}
where the tortuosity, $\tau$, is a function of the porosity in the Bruggemann approximation \cite{tjaden2016origin}
\begin{equation}\label{eq:bruggemann}
    \tau = \phi^{-0.5},
\end{equation}
and, as noted by Carman \cite{carman1937fluid} for spherical particles, $S_v = 6/D$.

The non-solid-phase comprises both liquid and gas phases with saturations $S_l$ and $S_g$, respectively.
These saturations reflect the respective volume fraction in the pores around the solid particles, with $S_l + S_g = 1$.
One notable simplification from the models considered by Voskuilen et al. \cite{voskuilen2021multi} is the use of constitutive models describing the liquid-phase saturation as a function of the liquid-phase pressure, $p_l$.
As a result, the gas phase staturation is calculated merely as $S_g = 1 - S_l$, allowing for the liquid-phase continuity to be decoupled from the gas phase continuity.
The implication of this is that the contribution from the gas phase pressure is omitted from consideration in the block multiphysics preconditioner in \Cref{sec:block-preconditioning}.

The liquid-phase obeys a continuity law of the form
\begin{equation}\label{eq:liquid-phase-continuity}
    \dfrac{\partial\left(c \phi S_l\right)}{\partial t} + \nabla \cdot\left(c \overrightarrow{v_l}\right) = S_{Li^+},
\end{equation}
where the liquid-phase velocity $\overrightarrow{v_l}$ is
\begin{equation}\label{eq:liquid-velocity}
    \overrightarrow{v_l} = -\dfrac{\kappa_i}{\mu_l}\left(\nabla p_l - C_m \nabla p_{l,0}\right),
\end{equation}
$C_m$ is a constant meant to prohibit liquid-phase transport below the melting temperature ($T_m$),
\begin{equation}\label{eq:regularization-term}
    C_m = \dfrac 1 2 \left(1 + \tanh{\left(\dfrac{T_m-T}{\Delta T}\right)}\right),
\end{equation}
and $c$ is the molar concentration, which itself is a constitutive function of temperature.
The continuity equation in \Cref{eq:liquid-phase-continuity} is solved for the liquid-phase pressure.
Due to the effects of capillary forces, the initial pressure gradient, $\nabla p_{l,0}$, is non-zero on the boundaries of materials.
As a result, Darcy's law is modified in \Cref{eq:liquid-velocity} that includes a regularization term of $C_m \nabla p_{l,0}$ to prevent advection before the electrolyte has melted \cite{voskuilen2016multi}.
Below the melting temperature, the liquid-phase viscosity, $\mu_l$, is further set to an arbitrarily large value to prevent advection.
Similar to the variable coefficient Poisson equation describing the solid-phase voltage in \Cref{eq:solid-phase-voltage},
the solution to the liquid-phase pressure, $p_l$, in \Cref{eq:liquid-phase-continuity,eq:liquid-velocity} can be re-arranged into a variable coefficient Poisson equation.
As a result, the same AMG techniques may be applied in solving for the liquid-phase pressure.

Within the liquid-phase, the transport equations for the mole fractions of $Li^+$ and $Cl^-$ are solved.
The mole fraction for $K^{+}$ is evaluated as $x_{K^+} = 1 - x_{Cl^-} - x_{Li^+}$
The transport equation for $x_{Li^+}$ resembles the liquid-phase continuity from \Cref{eq:liquid-phase-continuity},
but with a diffusive flux that is evaluated using the multi-component Stefan-Maxwell diffusion equation \cite{krishna1997maxwell} that depends on the gradients of the concentration and voltage,
\begin{equation}\label{eq:stefan-maxwell}
    \dfrac{\partial\left(c\phi S_l x_{Li^+}\right)}{\partial t} + \nabla\cdot \overrightarrow{D}_{Li^+} + \nabla\cdot\left(cx_{Li^+} \overrightarrow{v_l}\right)=S_{Li^+},
\end{equation}
This equation for the liquid-phase ion transport resembles an advection-diffusion equation, for which Petrov-Galerkin (PG) AMG may be suitable.
In practice, however, we observe that the domain-decomposition techniques discussed in \Cref{sec:domain-decomp} provide sufficiently robust preconditioners that exhibit little dependence on both the amount of compute resources (subdomains) or the problem scales considered in \Cref{chp:results}.

Finally, the concentration ($c_{s,i}$) for each solid-phase species $i$ is given by the ordinary differential equation (ODE)
\begin{equation}\label{eq:solid-phase-species-conduction}
    \dfrac{\partial c_{s,i}}{\partial t} = S_{s,i},
\end{equation}
where $S_{s,i}$ is a source term.
As a result, the subset of the monolithic linear system associated with the solid-phase species conduction is diagonal.
This fact is later exploited in \Cref{sec:block-preconditioning}, where the equations associated with this block are solved via a single sweep of Jacobi.

\mychp{Multiphysics Preconditioning}
\label{chp:preconditioning}
The performance of multiphysics simulations in \aria{} relies heavily on effective multiphysics preconditioning strategies.
We present two common approaches that are available within the \aria{}: `block-box' (physics agnostic) domain-decomposition preconditioners
and physics-aware block-based preconditioning.

Domain-decomposition approaches to multiphysics preconditioning are discussed in \Cref{sec:domain-decomp}.
Domain-decomposition techniques, while agnostic to the specific multiphysics use-case, remain attractive for preconditioning due to their user-friendly construction and excellent parallel scalability in applying local subdomain solvers.
At the same time, the inherently local construction of the one-level method means that the preconditioner may be ineffective for long, coarse modes that span several subdomains or the entire simulation domain.
As a consequence, the convergence rate of the domain-decomposition approaches rely heavily on the underlying multiphysics system and the number of processors used.

Block-based preconditioning strategies, in contrast, require comprehensive knowledge of the underlying physics.
In particular, a user must choose the block preconditioning strategy, as well as effective preconditioners or solvers for each of the subblocks that need to be inverted.
We propose a hierarchical block Gauss-Seidel preconditioner specially formulated for preconditioning the electrochemical solve phase for thermal battery simulations in \Cref{sec:block-preconditioning}.

\mysec{Domain-Decomposition Preconditioning}
\label{sec:domain-decomp}

In \aria{}, common choices for preconditioning multiphysics problems are domain-decomposition methods.
These methods are popular in part due to their effectiveness for a wide range of multiphysics applications and their user-friendly interface.
Due to their purely algebraic construction, the domain-decomposition methods considered are `black-box' approaches to multiphysics problems.
The \ifpack{} package \cite{prokopenko2016ifpack2} from \trilinos{} is used to configure a one-level domain-decomposition approach.
Subdomains are defined as processor-local submatrices, or small extensions thereof,
\begin{equation}\label{eq:processor-local-matrices}
   A_i = R_i A R_i^T,
\end{equation}
where $R_i$ is a boolean restriction matrix mapping degrees of freedom from the global system to the subdomain \cite{prokopenko2016ifpack2}.
The amount of overlap is left as a user-defined parameter, with no overlap being permissable.
An increase in the degree of overlap improves the convergence rate of the preconditioner \cite{toselli2004domain,dryja1990some},
albeit at the expense of additional communication and larger local subdomain problems.
Following Cai and Sarkis \cite{cai1999restricted}, we employ a restricted additive Schwarz (RAS) method as a preconditioner
\begin{equation}\label{eq:one-level-schwarz-exact-subdomain-solve}
   \tilde{M}_{RAS}^{-1} = \sum_{i=1}^{N} R_i^T D_i A_i^{-1} R_i,
\end{equation}
where $D_i$ is a weight matrix that is unity on the interior of the subdomain and zero in the overlap region.
The effect of this is that the subdomain solutions in the overlapping region are ignored, rather than added together.
As a consequence, the application of $R_i^T D_i$ avoids communication along the shared interface and/or overlapping regions and converges faster than an additive Schwarz scheme \cite{efstathiou2003restricted}.
The action of subdomain inverses is approximated using $M_i^{-1}$, an incomplete LU (ILU) factorization preconditioner on each subdomain.
The resulting RAS preconditioner is thus given by
\begin{equation}\label{eq:one-level-schwarz}
   M_{RAS}^{-1} = \sum_{i=1}^{N} R_i^T D_i M_i^{-1} R_i.
\end{equation}
The \emph{level-of-fill} \cite{saad2003iterative} is left as an option to the user, with higher levels of fill providing more accurate approximations to the action of $A_i^{-1}$.

While the one-level RAS preconditioner in \Cref{eq:one-level-schwarz} is effective in many multiphysics applications,
it nevertheless \emph{lacks} a coarse-grid correction, which is required to ensure scalability of the preconditioner \cite{toselli2004domain}.
Further, for certain multiphysics applications, one-level domain decomposition techniques are not effective.
For example, using of RAS as a preconditioner for the electrochemical solves considered in this report results in a lack of convergence, even at relatively small problem scales.
This occurs despite efforts to increase both the ILU fill level and subdomain overlap level.
In this context, even relying on a direct solver, such as those provided by \klu{}, are not sufficient to ensure convergence in a reasonable number of iterations.
Due to this lack of convergence, we omit any scaling studies utilizing the RAS preconditioner for the monolithic scaling studies in \Cref{chp:results}.
We will, however, see the use of RAS as a key component in conjunction with the preconditioner described in \Cref{sec:block-preconditioning}.

Several strategies could have been used to help mitigate the issues identified in the preceding paragraph.
For example, the introduction of an algebraic coarse grid correction provided through software packages such as \frosch{} \cite{heinlein2018frosch,heinlein2022parallel}.
Through the introduction of a reduced dimension GDSW coarse-grid \cite{dohrmann2017design}, Heinlein, Perego, and Rajamanickam demonstrated they could achieve iteration counts with little variation with respect to the number of subdomains \cite{heinlein2022frosch}.
In lieu of the approaches provided through the \frosch{} software package, we instead consider the construction of block-based, physics-aware preconditioning in \Cref{sec:block-preconditioning}.

\mysec{Block-based Preconditioning}
\label{sec:block-preconditioning}

The construction of block-based, physics-aware preconditioners requires the careful consideration of the various physics in the monolithic system and specific strategies for inverting the individual physics subblocks.
Because of their complexity, users of \aria{} often prefer the more user-friendly domain-decomposition approaches described in \Cref{sec:domain-decomp}.
In contrast to this approach, we consider a physics-based block preconditioner for thermal batteries.
Physics-based preconditioning strategies for solving electrochemical systems has been previously considered.
Similar approaches that employ block Gauss-Seidel based preconditioners have been considered by Fang et al. \cite{fang2019parallel} and Allen et al. \cite{allen2021segregated} in the context of modeling lithium-ion batteries.
The effort to construct a block-based preconditioner for the thermal battery applications is especially justified given the desire to provide robust, scalable solvers for thermal battery applications, wherein the solution to the electrochemical solve is the dominant bottleneck in an \aria{} simulation.

We employ the \teko{} package \cite{cyr2016teko} from \trilinos{} to construct a hierarchical block Gauss-Seidel preconditioner.
The full electrochemical system is described by a nested 2x2 block structure,
\begin{equation}\label{eq:block-splitting}
 \A = \begin{bmatrix}
    \A_{v,v} & \A_{v,n} \\
    \A_{n,v} & \A_{n,n}
 \end{bmatrix},
\end{equation}
where $\A_{n,n}$ represents the non-voltage matrix and $\A_{v,v}$ represents the voltage matrix.
The former system represents the discrete contribution from the Stefan-Maxwell diffusion equation terms for the liquid-phase species fractions from \Cref{eq:stefan-maxwell},
the liquid-phase continuity from \Cref{eq:liquid-phase-continuity}, and the solid-phase species conduction terms from \Cref{eq:solid-phase-species-conduction}.
The voltage matrix $\A_{v,v}$, on the other hand, only comprised two degrees-of-freedom: the liquid-phase voltage from \Cref{eq:liquid-phase-voltage} and the solid-phase voltage from \Cref{eq:solid-phase-voltage}.
The cross-coupling terms between the non-voltage and voltage blocks are in $\A_{v,n}$ and $\A_{n,v}$.
The coupled solid/liquid-phase voltage system is given by
\begin{equation}\label{eq:voltage-block}
 \A_{v,v} = \begin{bmatrix}
    A_{\Phi_s, \Phi_s} & A_{\Phi_s, \Phi_l} \\
    A_{\Phi_l, \Phi_s} & A_{\Phi_l, \Phi_l}
 \end{bmatrix},
\end{equation}
where $A_{\Phi_s, \Phi_s}$ and $A_{\Phi_l, \Phi_l}$ represent the solid and liquid-phase voltages, respectively.
$A_{\Phi_s, \Phi_l}$ and $A_{\Phi_l, \Phi_s}$, on the other hand, represent the coupling of the two phases through the Butler-Volmer potential form in \Cref{eq:butler-volmer}.
The non-voltage terms are described by
\begin{equation}\label{eq:non-voltage}
 \A_{n,n} = \begin{bmatrix}
    A_{s} & & \\
    & A_{x,x} & A_{x,p} \\
    & A_{p,x} & A_{p,p}
 \end{bmatrix},
\end{equation}
where $A_{s}$ is a diagonal matrix corresponding to the solid-phase species conduction (uncoupled ODEs) from \cref{eq:solid-phase-species-conduction},
$A_{x,x}$ and $A_{p,p}$ correspond to the liquid-phase species fractions and liquid-phase pressure, respectively,
and $A_{x,p}$ and $A_{p,x}$ are the coupling blocks between them.

Since the coupled voltage terms and non-voltage terms are somewhat weakly coupled,
we consider the following hierarchical block Gauss-Seidel preconditioner
\begin{equation}\label{eq:hierarchical-bgs-exact-inverse}
 \M^{-1} = \begin{bmatrix}
    \A_{v,v} & \A_{v,n} \\
    & \A_{n,n}
 \end{bmatrix}^{-1}.
\end{equation}
Provided we employ a sparse direct solver for the required action of $\A_{n,n}^{-1}$ and $\A_{v,v}^{-1}$ in the preconditioner,
convergence is reached in very few iterations.
This approach, however, is not scalable.
We instead employ GMRES(30) to invert the $\A_{n,n}$ and $\A_{v,v}$ blocks using $\M_{n,n}^{-1}$ and $\M_{v,v}^{-1}$ as preconditioners, respectively.
A $10^{-6}$ relative residual tolerance for both blocks is used.
For the non-voltage block, we employ as a preconditioner a block Gauss-Seidel scheme
\begin{equation}\label{eq:non-voltage-bgs}
 \M_{n,n}^{-1} = \begin{bmatrix}
    A_{s} & & \\
    & M_{x,x} & A_{x,p} \\
    & & M_{p,p}
 \end{bmatrix}^{-1}.
\end{equation}
$A_{s}$ is a diagonal matrix, so we directly invert it.
$M_{x,x}^{-1}$ and $M_{p,p}^{-1}$ are preconditioners for $A_{x,x}$ and $A_{p,p}$, respectively.
For the liquid-phase species fractions, the subblock is preconditioned via the DD(0)-ILU(0) scheme from \ifpack{} as described in \Cref{sec:domain-decomp}.
For the liquid-phase pressure, we employ a \muelu{} smoothed aggregation algebraic multigrid (SA-AMG) V-cycle using a second-order Chebyshev smoother as a preconditioner.
More details regarding SA-AMG are deferred to \Cref{sec:sa-amg}.
Similar to the to the non-voltage block in \Cref{eq:non-voltage-bgs},
we employ a 2x2 block Gauss-Seidel scheme for the coupled voltage equation
\begin{equation}\label{eq:voltage-bgs}
 \M_{v,v}^{-1} = \begin{bmatrix}
    M_{\Phi_s, \Phi_s} & A_{\Phi_s, \Phi_l} \\
     & M_{\Phi_l, \Phi_l}
 \end{bmatrix}^{-1}.
\end{equation}
For the solid and liquid voltage preconditioners $M_{\Phi_s, \Phi_s}^{-1}$ and $M_{\Phi_l, \Phi_l}^{-1}$, we employ SA-AMG.
Despite the exponential coupling term between the solid and liquid-phase voltages in the Butler-Volmer potential in \Cref{eq:butler-volmer},
the two-phase voltage system in \Cref{eq:voltage-block} does not require a monolithic AMG approach, such as those considered by Ohm et al. in \cite{ohm2021monolithic} or Fang et al. in \cite{fang2019parallel}.

The two inner GMRES solves invoked in each preconditioner evaluation in \Cref{eq:hierarchical-bgs-exact-inverse} require the use of \emph{flexible} GMRES \cite{saad1993flexible} to solve \Cref{eq:block-splitting}.
The only cost increase associated with flexible GMRES over GMRES is doubling the memory requirement associated with storing the restart space.
However, as previously noted, the preconditioner described in \Cref{eq:hierarchical-bgs-exact-inverse} converges in few iterations.
For example, flexible GMRES(5) suffices.
Additionally, flexible GMRES actually requires one fewer preconditioner evaluation by omitting the preconditioner evaluation required in forming the approximate solution after the Arnoldi process \cite{saad1993flexible}.
For a high-quality preconditioner that converges in very few ($<5$) iterations, this represents a significant cost savings.

\mysubsec{Smoothed Aggregation Algebraic Multigrid}
\label{sec:sa-amg}

The liquid-phase continuity, liquid-phase voltage, and solid-phase voltage subblocks described in the block-based multiphysics preconditioner require the use of algebraic multigrid (AMG) techniques.
All three of these blocks share similar AMG configurations.
We therefore conclude our discussion on block-based multiphysis for thermal battery applications by describing the smoothed aggregation (SA) AMG methods needed.
Only a very high-level description of SA-AMG is presented here.
For AMG, see \cite{BrMcRu84,RuSt85,Falgout2006a} and references therein.
For SA, see \cite{Vanek1996a,Vanek2001a}.

Algebraic multigrid (AMG) methods \cite{BrMcRu84,RuSt85,Falgout2006a} are a popular and effective class of methods for solving systems arising from elliptic PDEs.
Key components of AMG methods include prolongators, restrictors, coarse-level operators, and smoothers.
The coarse level operators are constructed in the usual Galerkin formulation with
\begin{equation}\label{eq:galerkin}
   A_0 = A,\ A_{l+1} = R_l A_l P_l, l=1,\dots,\ell-1,
\end{equation}
where $P_l$ is the prolongation operator mapping degrees of freedom from the coarser grid to the fine grid,
$R_l=P_l^T$ is the restriction operator mapping degrees of freedom from the fine grid to the coarser grid,
and $A_{l+1}$ is the coarser level operator.
\emph{Smoothed aggregation} AMG, moreover, seeks to minimize the energy of the coarse basis functions by smoothing the tentative prolongation operator to construct $P_l$ \cite{Vanek1996a,Vanek2001a}.

We use SA-AMG as provided by the \muelu{} package \cite{BergerVergiat2023a} in \trilinos{}.
Chebyshev polynomial smoothing is an effective parallel multigrid smoother \cite{adams2003parallel,baker2011multigrid} that can even be used to accelerate the one-level Schwarz methods described in \Cref{sec:domain-decomp} \cite{phillips2022tuning,phillips2025optimal}.
While generally tailored for symmetric positive definite (SPD) systems, the Chebyshev smoothing method can be adapted for non-symmetric problems by forming an ellipse to bound the convex hull of the spectra \cite{manteuffel1977tchebychev}.
Unless otherwise noted, we use a 2nd-order Chebyshev polynomial smoother, which corresponds to two sweeps as pre- and post smoother on each level.
The coarsest level is solved with a \klu{} direct method as provided by the \amesos{} package in \trilinos{} \cite{Bavier2012a}.
In order to control the operator complexity, matrix filtering is employed during smoothing of the tentative prolongator \cite{Gee2009a},
with a diagonal weight matrix based on the absolute row sum used to better capture the scaling of each row in the filtered matrix \cite{Hu2022a} .
The coordinate based multijagged approach in \zoltan{} \cite{Deveci2015} is used to rebalance coarser level ($l > 1$) matrices, targeting \num{10000} rows per processor with no fewer than \num{1000} rows per processor.
Levels are generated until reaching a target maximum coarse grid size of \num{5000}.

\mychp{Results}
\label{chp:results}
We evaluate the strong and weak scalability of the preconditioners described in \Cref{chp:preconditioning} for a challenging thermal battery configuration.
All results utilize the Eclipse system at Sandia National Laboratories, which features a dual socket, 18 core Intel Broadwell E5-2695 processor per node.
This corresponds to 36 MPI ranks per node.
Even at the smallest problem scales considered, the `black-box' DD-ILU based preconditioners described in \Cref{sec:domain-decomp} fail to converge on the electrochemical system in \cref{eq:block-splitting}.
As such, we focus our attention on the block-based multiphysics preconditioner described in \Cref{sec:block-preconditioning}.
We consider the solution to the electrochemical system in \cref{eq:block-splitting} for a fixed snapshot taken during the electrochemical activation for the 2D axisymmetric simulation of thermal batteries for four different mesh resolutions.
The initial coarse mesh considered in the TABS repository \cite{voskuilen2021multi,roberts2023development,roberts2018tabs} represents a coarse resolution featuring multiple elements along the axial dimension to resolve the multiple interfaces between the cathode, separator, and anode layers within the battery.
In the radial dimension, however, very few elements are used.
As a consequence, the element aspect ratios are severely stretched, with the aspect ratio varying from 11.28 to 37.03.
The anisotropic mesh geometry of this case combined with the relatively small problem size makes this a challenging case to start evaluating solver scalability.
In this work, we opt for a 16-fold refinement in the number of elements in the radial dimension while using the same number of elements in the axial dimension.
In this case, the element aspect ratios are improved, varying from 1.17 to 2.78.
This forms our baseline ($N=1$) problem scale with roughly $E=126.9K$ elements and $n=843.6K$ degrees of freedom in the monolithic system.
From this point, three larger cases are generated by performing succesive levels of uniform mesh refinment (UMR), bringing the largest case ($N=64$) to $E=8.1M$ elements and $n=51.3M$ degrees of freedom.

Strong scaling results are generated by considering the performance of the baseline scale with $E=126.9K$ elements on a single processor.
From here, successive processor doublings are done until reaching $P=2^6$ processors.
The same is repeated for each problem scale $N$, wherein $P=N$ ranks are initially used until reaching $P=2^6\ N$ ranks.
The scaling study thus spans between a single processor to $P=2^{12}$ ranks.
This simulatenously allowes weak scaling performance to be gleaned from lines of constant $E/P$.

The results are separated into two sections: \Cref{sec:subblock-solvers} and \Cref{sec:end-to-end-solver}.
The former looks at the scalability of the various preconditioners required in inverting the subblocks for coupled non-voltage preconditioner in \cref{eq:non-voltage-bgs}.
In addition, scalability of the individual voltage equations is assessed as a stepping stone to coupled voltage solve in \Cref{eq:voltage-block}.
The latter section combines then looks at the scalability and convergence for the entire end-to-end monolithic solve.

We note the similarity of the scaling plots in \Cref{fig:liquid-mass-fraction,fig:liquid-phase-pressure,fig:liquid-voltage,fig:solid-voltage,fig:coupled-voltage,fig:end-to-end-solve}, so we describe them here.
The setup and solve phases of each solver are repeated 100 times for each configuration, with the first run discarded as an initial warmup.
The mean setup, solve, and combined setup plus solve times are reported, with the error bars representing one standard deviation from the mean.

Each color in \Cref{fig:liquid-mass-fraction,fig:liquid-phase-pressure,fig:liquid-voltage,fig:solid-voltage,fig:coupled-voltage,fig:end-to-end-solve} corresponds to a different level of mesh refinement: blue corresponds to the 16x radial refinement ($N=1$),
orange represents a single UMR of the case in blue ($N=4$),
green represents a single UMR of the case in orange ($N=16$),
and red represents a single UMR of the case in green ($N=64$).
For each problem scale, we consider running at several element-to-processor ratios, $E/P$.
At fixed $E/P$ ratios, moving from the unrefined case represented by the blue curve to the refined case in red represents \emph{weak} scaling.
We use as our weak scaling model
\begin{equation}\label{eq:weak-scaling-efficiency}
    T_n = \dfrac{T_{m}}{\etaweak{}^{\log_2\left(\dfrac{n}{m}\right)}},
\end{equation}
where $T_n$ denotes the execution time for a problem with $n$ degrees of freedom, $T_m$ denotes the execution time for $m$ degrees of freedom with $m < n$, and $\eta$ represents the weak scaling parallel efficiency.
The model in \cref{eq:weak-scaling-efficiency} requires a \emph{fixed} number of degrees of freedom per rank, i.e., $n/P$ is constant.
Perfect weak scaling ($\etaweak{} = 1$) corresponds to a flat line, i.e., setup and solution time have no dependence on problem scale.
An efficiency of $\etaweak{}=0.5$, on the other hand, indicates that $T_n = 2 T_m$ when $n = 2m$.
That is, for every doubling in problem size, the time also doubles.
A gray dashdot curve representing the weak scaling efficiency across the largest three problem scales at the largest $N/P=4$ ratio is shown with the scaling efficiency noted.
A least squares regression from Scipy \cite{virtanen2020scipy} is used to estimate the efficiency based on the weak scaling cost model described in \cref{eq:weak-scaling-efficiency}.
The reasons for this particular reporting of the weak scaling efficiency are three fold:
first, the communication stencil for the matrix-vector product of a two dimensional $\mathbb{Q}_1$ finite element Poisson problem,
which is required in the Krylov solver and smoother and residual evaluations for a multigrid preconditioner,
requires a minimum of 9 processors to fully saturate;
second, the processor count required to reach this saturation point for the baseline $N=1$ problem scale would surpass the strong scale limit for most of the subblock solvers considered in \Cref{sec:subblock-solvers};
and finally, since each node of Eclipse only has 36 MPI ranks, the change in problem scale necessitates moving from purely on-node communication to off-node communication;
The weak scaling at this problem scale therefore better illustrates the real conditions under which a user would percieve the parallel scalability of the solver.

The problem scales described in the preceding paragraph are also used to assess the strong scalability of the problem.
Along lines of fixed problem scale (color in \Cref{fig:liquid-mass-fraction,fig:liquid-phase-pressure,fig:liquid-voltage,fig:solid-voltage,fig:coupled-voltage,fig:end-to-end-solve}), \emph{strong} scaling is assessed as the total number of ranks $P$ increases.
We use as our strong scaling model
\begin{equation}\label{eq:strong-scaling-efficiency}
    T_{P_n}= \dfrac{1}{\etastrong{}}\left(\dfrac{P_m}{P_n}\right)T_{P_m},
\end{equation}
where $T_{P_n}$ denotes the time using $P_n$ processors, $T_{P_m}$ is the time required using a baseline configuration of $P_m$ processors, and $\etastrong{}$ is the strong scaling efficiency.
Ideal strong scaling ($\etastrong{}=1$) is represented with an idealized gray dashed curve.
In this case, doubling the processor count halves the execution time.
For suboptimal strong scaling with $\etastrong{} = 0.5$, a doubling in the processor count reduces the execution time by only 25\%.
Fischer et al. characterize the configuration at which further strong scaling reaches a point of diminishing returns as defined by $\etastrong{} = 0.8$ \cite{fischer2020scalability}.
We employ a similar metric here, but use $\etastrong{} = 0.5$ as our cut-off for strong scaling.

\mysec{Scalable Subblock Solvers}\label{sec:subblock-solvers}

In this section, we assess the convergence and parallel scalability of the subblock preconditioners described in \Cref{sec:block-preconditioning}.
This includes the solution to the following:
\begin{itemize}
\item \Cref{sec:liquid-mass-fraction}: liquid-phase species fractions, which are described via the multicomponent Stefan-Maxwell diffusion model in \Cref{eq:stefan-maxwell},
\item \Cref{sec:liquid-phase-continuity}: liquid-phase pressure, which is described via the liquid-phase continuity equation in \Cref{eq:liquid-phase-continuity},
\item \Cref{sec:liquid-phase-voltage}: liquid-phase voltage, which is described with the diffusion-dominated multi-component advection-diffusion equation in \Cref{eq:liquid-phase-voltage},
\item \Cref{sec:solid-phase-voltage}: solid-phase voltage, which is described with the Poisson equation in \Cref{eq:solid-phase-voltage}, and
\item \Cref{sec:coupled-voltage}: coupled voltage system, which is both the solid and liquid-phase voltages plus coupling terms in \Cref{eq:voltage-block}.
\end{itemize}
Each system noted above is solved to a relative residual criteria of $\dfrac{\norm{\underline r^N}}{\norm{\underline r^0}} \leq 10^{-8}$, unless otherwise noted.
Iteration counts for each subblock solver and the end-to-end electrochemical solve at the strong scaling limit ($N/P=8$) are summarized in \Cref{tab:iteration-counts}.
In practice, only the coupled voltage system and coupled non-voltage system are solved during the preconditioner evaluation in \Cref{eq:hierarchical-bgs-exact-inverse}.
The preconditioners for the liquid-phase species fractions and liquid-phase pressure are employed in a block-based preconditioner for the non-voltage terms in \Cref{eq:non-voltage-bgs}.
Further, the multigrid preconditioners for the liquid and solid-phase voltages are used only as an intermediate step to construct the block preconditioner described in \Cref{eq:voltage-bgs}.
Nevertheless, this demonstrates that appropriate preconditioners are used for each subblock, both in terms of convergence and parallel scalability.

Finally, we note that we have omitted any scalability studies for solving the solid-phase species conduction subblock in \Cref{eq:solid-phase-species-conduction}.
As previously mentioned, this system is diagonal, and therefore the $A_s^{-1}$ evaluation required by the non-voltage block preconditioner in \Cref{eq:non-voltage-bgs} requires only a simple diagonal scaling to solve.

\mysubsec{Stefan-Maxwell Diffusion Liquid-Phase Species Fractions}\label{sec:liquid-mass-fraction}

\begin{figure}[h]
    \includegraphics[width=\textwidth]{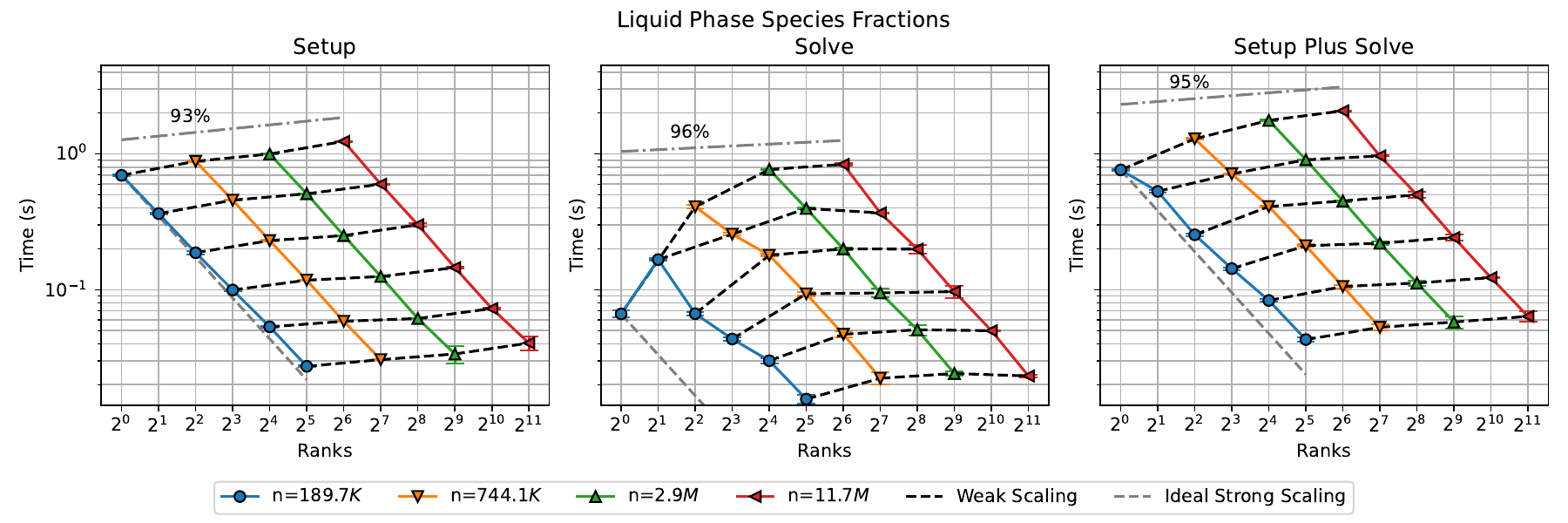}
    \caption{
        \label{fig:liquid-mass-fraction}
        Scalability of GMRES preconditioned via DD(0)-ILU(0) for the liquid-phase species fractions.
        Strong scalability is maintained up to $n/P\sim\num{6000}$, with sub-second setup and solve times achieved across all problem scales.
        The weak scaling efficiency reaches 95\% for the combined setup and solve phase, highlighting the efficiency of the preconditioning strategy for this subblock.
    }
\end{figure}

While the domain-decomposition preconditioners described in \Cref{sec:domain-decomp} are unable to solve the monolithic system,
they nevertheless prove effective as preconditioners for individual subblocks.
The liquid-phase species fractions in \cref{eq:stefan-maxwell} resemble an advection-diffusion problem,
which is amenable to Petrov-Galerkin AMG (PG-AMG) \cite{idelsohn1996petrov,wu2006analysis} or multilevel Schwarz \cite{lin2010parallel} preconditioning approaches.
As noted by Lasser and Toselli \cite{lasser2003overlapping}, the introduction of a coarse-grid space to a one-level Schwarz method for advection-diffusion problems improves the iteration count.
We employ here a single-level non-overlapping restricted additive Schwarz preconditioner as noted in \Cref{eq:one-level-schwarz} with a single sweep of ILU(0) from \ifpack{} as a subdomain solver.
A GMRES(30) solver from \belos is employed until an 8 order-of-magnitude reduction in the residual is achieved.
With this preconditioner, convergence is achieved in 27 or fewer iterations, regardless of the problem scale or number of processors (subdomains).
For $P>1$, the iteration count shows little dependence on problem scale or the number of subdomains.
The preconditioner setup cost, solve time, and combined setup and solve are reported in \Cref{fig:liquid-mass-fraction}.

With exception to the $n=189.7K$ line, \Cref{fig:liquid-mass-fraction} demonstrates that the setup, solve, and combined times have excellent strong scalability.
In the smallest scale case, a large increase in the solve phase is observed going from $P=1$ to $P=2$ ranks.
In this case, the preconditioner changes from a ILU(0) preconditioner to a domain-decomposition preconditioner, leading to an increase in the overall iteration count required to reach convergence.
For the other problem scales considered, the slope of the strong scaling curves closely resembles the ideal strong line.
This result is expected, provided the iteration count remains somewhat low and fixed with respect to increasing $P$ (and thus the number of subdomains).
The communication patterns required in the domain-decomposition preconditioner and matrix-vector product require only highly-scalable point-to-point communication with the values in neighboring elements.
That is, once the number of processors exceeds the point at which the communication stencil is fully saturated, overhead associated with the communication costs of these operators will scale as $\mathcal O(1)$ with respect to $P$.
At the same time, if the iteration count is somewhat large and grows with respect to $P$, the strong scalability will be poor.
This is due to the $\mathcal O(\log(P))$ communication complexity associated with the all-reduce \cite{thakur2005optimization,hoefler2010toward} required in the orthogonalization step in (flexible) GMRES.
No more than 27 iterations are ever needed to meet 8 order-of-magnitude reduction required in the residual.
The weak scalability of the solve and setup phases is similarly quite good, achieving 93\% and 96\% parallel efficiency, respectively.
This equates to a 95\% weak scaling parallel efficiency for the combined setup and solve phase.

The strong scalability of the setup, solve, and combined phases is similarly excellent.
For example, the strong scaling efficiency for the setup and solve phases at the $n=2.9M$ and $n=11.7M$ scales exceeds 92\% across all of the processor counts considered.
As such, the strong scaling limit ($\etastrong{} \leq 0.5$) has not been reached in \Cref{fig:liquid-mass-fraction}, despite the somewhat aggresive $n/P \sim \num{6000}$ used at the highest processor counts.
While additional compute resources can be effectively used to solve this particular subblock, we note that the other subblock solver components will reach their respective strong scale limits at much larger values of $n/P$.
We lastly note that the little time is required for the overall solution to this subblock.
Sub-second solves are achieved irrespective of the problem scale or resources used.
The total end-to-end setup plus solve time reaches sub-second for all problem scales, provided $n/P < \num{100000}$.
Finally, since this subblock strong scales so well, the end-to-end setup plus solve time is achieved in less than 0.1 seconds at the largest processor counts.

\mysubsec{Liquid-Phase Continuity}\label{sec:liquid-phase-continuity}

\begin{figure}[h]
    \includegraphics[width=\textwidth]{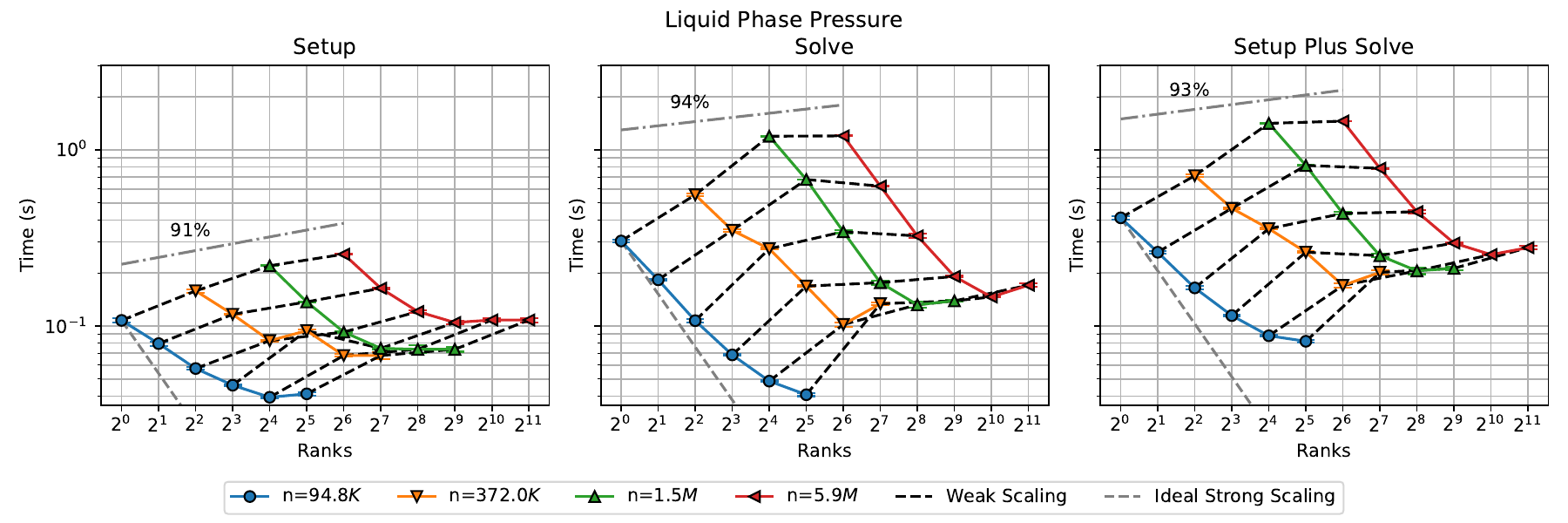}
    \caption{
        \label{fig:liquid-phase-pressure}
        Scalability of GMRES(30) preconditioned via SA-AMG for the liquid-phase pressure.
        Setup and solve times remain below 0.3 seconds at larger processor counts, with strong scalability observed up to $n/P\sim\num{11500}$.
        The weak scaling efficiency of 93\% demonstrates the robustness of the approach for handling the liquid-phase pressure equations.
    }
\end{figure}

The SA-AMG approach described in \Cref{sec:sa-amg} from \muelu{} is employed as a preconditioner for a GMRES(30) solver from \belos{}.
The results are reported in \Cref{fig:liquid-phase-pressure}.
The strong scalability characteristics differ from the results presented in \Cref{sec:liquid-mass-fraction}.
While the results in the preceding subsection were based on highly scalable processor-local subdomains with little communication required,
the algorithms employed in SA-AMG require more complex communication patterns,
such as the intrinstric all-to-all pattern required in the coarse grid solve \cite{bello2022polynomial}.
As a consequence, we observe limits in the strong scalability of this approach relative to the one-level Schwarz approach in the preceding subsection.
The setup phase reaches the $\etastrong{} \leq 0.5$ limit at $n/P\sim\num{23000}$ and shows no significant improvement beyond $n/P\sim\num{11500}$.
The solve phase, which dominates the overall cost, shows better strong scaling characteristics.
The strong scaling limit for the solve is reached with $n/P\sim\num{6000}$ with a solution time less than 0.2 seconds for all problem scales.
The combined setup plus solve thereby reaches the scaling limit at an intermediate $n/P\sim\num{11500}$.
At this point, the end-to-end solution time is less than 0.3 seconds.
While the subblock in the preceding subsection is larger and strong scales better, the end-to-end solution time achieves a similar time near 0.2 seconds at this processor count.
The weak scalability (\Cref{eq:weak-scaling-efficiency}) of the liquid-phase continuity setup and solve is excellent, with 91\% efficiency in the setup and 94\% efficiency in the solve phase.
Combined, this yields a 93\% weak scaling efficiency for the end-to-end setup plus solve.

\mysubsec{Liquid-Phase Voltage}\label{sec:liquid-phase-voltage}

\begin{figure}[h]
    \includegraphics[width=\textwidth]{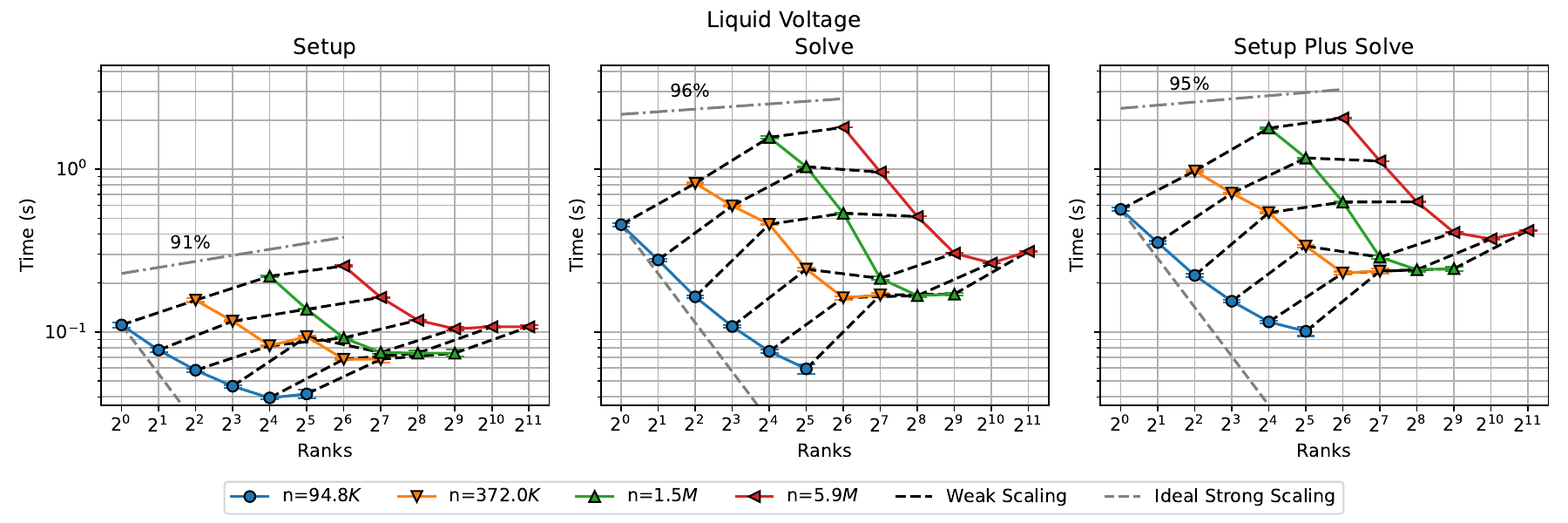}
    \caption{
        \label{fig:liquid-voltage}
        Scalability of GMRES(30) preconditioned via SA-AMG for the liquid-phase voltage.
        End-to-end setup and solve times remain below 0.4 seconds, with strong scalability achieved up to $n/P\sim{11500}$.
        The weak scaling efficiency of 95\% confirms the effectiveness of the approach for solving diffusion-dominated equations in this subblock.
    }
\end{figure}

\Cref{eq:liquid-phase-voltage}, the liquid-phase voltage equation, does not resemble a Poisson equation at initial glance.
However, the advective $\tau\dfrac{\partial\Phi_l}{\partial t}$ term is small thanks to the small magnitude of $\tau$.
This indicates that the equation is similar to a diffusion-dominated advection-diffusion equation due to the diffusive $\nabla\cdot\left(F\sum_j \overrightarrow{D}_j z_j\right)$ term.
In fact, the diffusion term is so dominant in this case that SA-AMG is sufficient for the solution of this subblock without needing to resort to a PG-AMG approach.
Much like the configuration for the liquid-phase continuity solve in \Cref{sec:liquid-phase-continuity}, the SA-AMG approach described in \Cref{sec:sa-amg} from \muelu{} is employed as a preconditioner for a GMRES(30) solver from \belos{}.
Results are reported in \Cref{fig:liquid-voltage}.

The strong scalability of the solver closely mirrors that of the previous subsection.
The strong scalability of the setup phase hits the $\etastrong{} \leq 0.5$ limit at $n/P\sim\num{23000}$.
Past $n/P\sim\num{11500}$, the setup phase shows no benefit with respect to using additional compute resources.
As before, the solve phase is the majority of the overall end-to-end cost, but shows better strong scaling than the setup phase.
The strong scaling limit is reached with $n/P\sim\num{6000}$, leading to a solution time less than 0.3 seconds for all problem scales.
For the combined setup plus solve, the strong scaling limit occurs at $n/P\sim\num{11500}$.
At this point, the end-to-end solution time is less than 0.4 seconds.
The weak scalability of the liquid voltage setup and solve is excellent, with 91\% efficiency in the setup and 96\% efficiency in the solve phase
for a combined 95\% efficiency in the end-to-end setup plus solve.
Overall, the liquid-phase voltage subblock behaves comparably to the liquid-phase continuity from the preceding subsection.

\mysubsec{Solid-Phase Voltage}\label{sec:solid-phase-voltage}

\begin{figure}[h]
    \includegraphics[width=\textwidth]{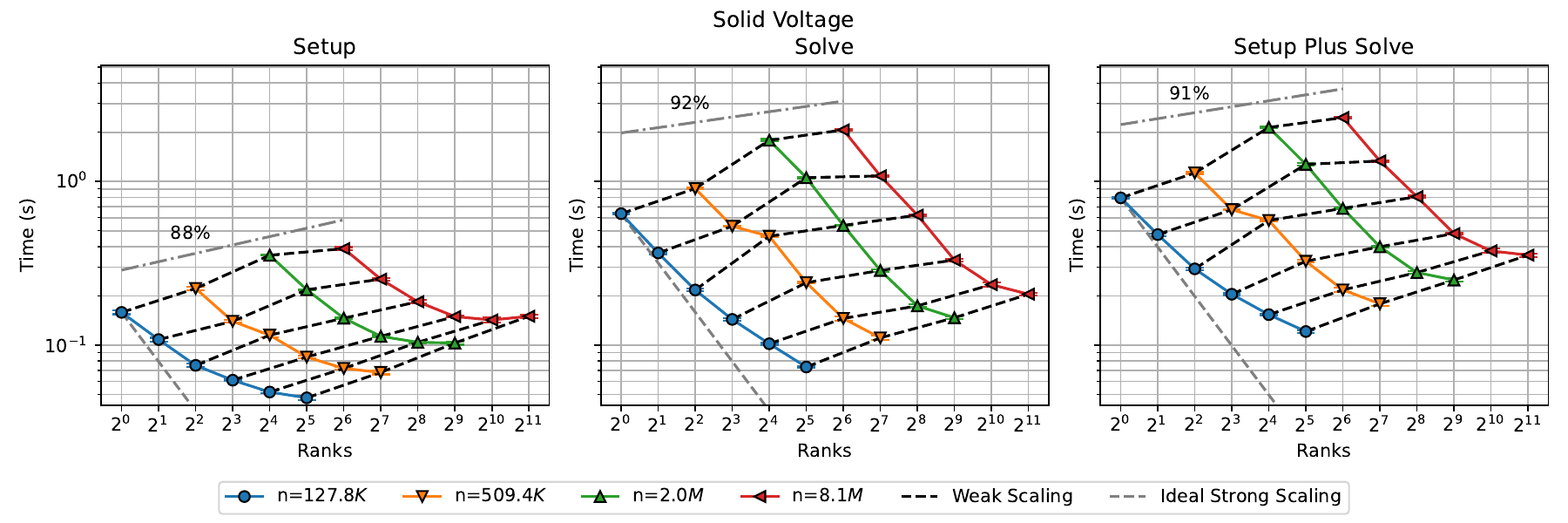}
    \caption{
        \label{fig:solid-voltage}
        Scalability of GMRES preconditioned via SA-AMG for the solid-phase voltage.
        Strong scalability is observed up to $n/P\sim\num{16000}$, with a combined setup plus solve time below 0.48 seconds.
        The weak scaling efficiency of 91\% highlights the ability of the approach to handle the challenges posed by highly variable electrical conductivity in the solid-phase voltage equation.
    }
\end{figure}

\Cref{eq:solid-phase-voltage}, which describes the solid-phase voltage equation, represents a variable coefficient Poisson equation, making this subblock especially amenable to SA-AMG approaches.
GMRES(30) preconditioned via SA-AMG is employed as the solver in this case, similar to the liquid-phase continuity in \Cref{sec:liquid-phase-continuity} or liquid-phase voltage in \Cref{sec:liquid-phase-voltage}.
However, the electrical conductivity, $\sigma$, varies by as much as 10 orders-of-magnitude across the simulation domain.
For example, $\sigma \approx 10^{-4}$ \unit{\siemens\per\meter} in the separator, while the conductivity is as large as $\sigma \approx 10^{6}$ \unit{\siemens\per\meter} in the anode.
As such, special attention is paid to to tuning the drop tolerance, $\theta$, used by the SA.
A grid search between $\theta = 0$ (no dropping) to $\theta = 0.64$ identified $\theta = 0.04$ as the optimal value, minimizing iteration count and cost for the baseline $N=1$ case.
Notably, this led to over a three-fold reduction in both iteration count and cost over the choice of no dropping with $\theta = 0$.

The setup phase exhibits similar poor strong scaling seen in the two preceding subsection, with $\etastrong{} \leq 0.5$ at $n/P\sim\num{32000}$.
The setup time at the largest problem scale ($N=64$) is 0.18 seconds.
The solve phase, in contrast, reaches its strong scaling limit at $n/P\sim\num{8000}$.
This corresponds to a solve time of 0.23 seconds or less.
Despite the better strong scaling, the time associated with the solve phase remains the majority of the end-to-end cost.
For the end-to-end setup plus solve time, the scaling limit is reached at $n/P\sim\num{16000}$, taking no more than 0.48 seconds at the largest problem scale.
For reference, this is only slightly more expensive than the roughly 0.4 seconds required for the liquid-phase voltage in the preceding subsection, despite the fact that the problem size considered here contains over 2 million more rows and is using the same processor count.
The weak scalability of the solid voltage setup and solve is similarly quite good, although somewhat less so than that for the liquid voltage.
Weak scaling efficiencies of 88\% and 92\% are reached in the setup and solve phases, respectively.
The end-to-end setup plus solve weak scaling efficiency is 91\%.

\mysubsec{Coupled Voltage System}\label{sec:coupled-voltage}

\begin{figure}[h]
    \includegraphics[width=\textwidth]{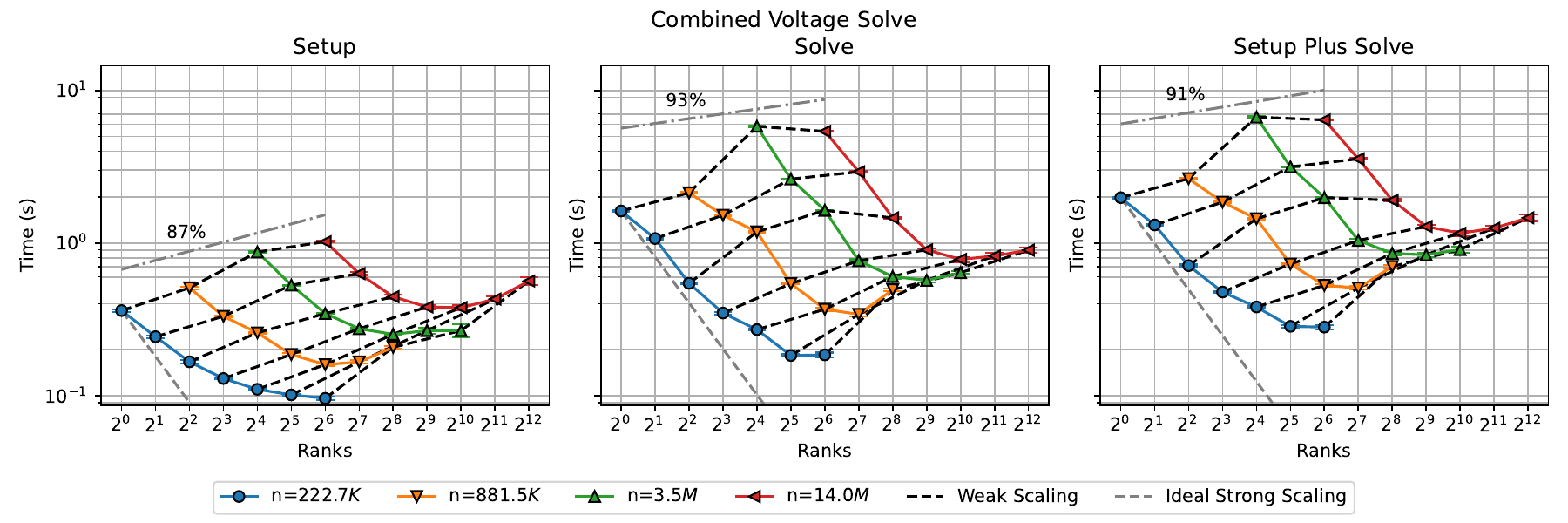}
    \caption{
        \label{fig:coupled-voltage}
        Scalability of the block Gauss-Seidel solver for the coupled solid/liquid voltage system.
        End-to-end setup and solve times remain below 1.28 seconds, with strong scalability achieved up to $n/P\sim{27500}$.
        The weak scaling efficiency of 91\% demonstrates the robustness of the preconditioner in managing the coupling between the solid and liquid-phase voltage equations.
    }
\end{figure}

With subblock solvers identified for the liquid and solid-phase voltages identified in \Cref{sec:liquid-phase-voltage,sec:solid-phase-voltage},
the coupled voltage system in \Cref{eq:voltage-block} may now be preconditioned using the block Gauss-Seidel scheme in \Cref{eq:voltage-bgs}.
A single sweep of the SA-AMG preconditioner from the solid voltages in \Cref{sec:solid-phase-voltage} and liquid voltages in \Cref{sec:liquid-phase-voltage}
are used as the action of $M_{\Phi_s, \Phi_s}^{-1}$ and $M_{\Phi_l, \Phi_l}^{-1}$ required in applying a block Gauss-Seidel sweep from \Cref{eq:voltage-bgs}.
One notable modification is made to the individual preconditioners in this section.
In lieu of a 2nd-order Chebyshev polynomial smoother from the two preceding sections, a 4th-order Chebyshev polynomial smoother is used to ensure that each individual subblock solver provides a sufficient reduction in the residual upon each sweep.
This corresponds to four sweeps as pre- and post smoother on each level, which makes $M_{\Phi_s, \Phi_s}^{-1}$ and $M_{\Phi_l, \Phi_l}^{-1}$ relatively heavy subblock preconditioners.
The block Gauss-Seidel preconditioner in \Cref{eq:voltage-bgs} is then used with a flexible GMRES(30) solver from \belos{} as the solver until a relative residual criterion of $\dfrac{\norm{\underline r^N}}{\norm{\underline r^0}} \leq 10^{-6}$ is reached.

We expect that this coupled voltage system solve will outweigh the non-voltage coupled solve required to apply one preconditioner sweep of the full, monolithic preconditioner in \Cref{eq:hierarchical-bgs-exact-inverse}.
There are two reasons for this.
First, the application of one sweep of \Cref{eq:voltage-bgs} requires two SA-AMG preconditioner evaluations that occur sequentially.
In contrast, the application of one sweep of the non-voltage block Gauss-Seidel preconditioner in \Cref{eq:non-voltage-bgs} requires only a single SA-AMG sweep for the liquid-phase continuity,
then a highly scalable DD(0)-ILU(0) sweep for the liquid-phase species fractions and a Jacobi sweep for the solid-phase species conduction.
The second reason for the relative difficulty of the coupled voltage system comes from the physics itself through the exponential interface coupling term between the voltage in the liquid and solid-phases in the Butler-Volmer potential in \Cref{eq:butler-volmer}.
Therefore, we defer discussion of the scalability of the coupled non-voltage solve to \Cref{sec:end-to-end-solver}.
Results for the coupled voltage system are shown in \Cref{fig:coupled-voltage}.

Strong scaling for the setup phase suffers from the same relative lack of scalability as the solvers noted in \Cref{sec:liquid-phase-voltage,sec:solid-phase-voltage},
with the $\etastrong{}\leq0.5$ scaling limit being hit as early as $n/P\sim\num{55000}$.
As expected, this number matches the sum of the strong scaling limits as $n/P\sim\num{32000}$ for the solid voltage and $n/P\sim\num{23000}$.
The most expensive setup cost occurs at the largest scale ($N=64$), where the setup takes 0.45 seconds.
While the bulk of the end-to-end setup plus solve cost is from the solve phase, its scalability is better than that of the setup step.
$\etastrong{}\leq0.5$ occurs past $n/P\sim\num{27500}$, where the solve phase takes as much as 0.9s.
Since the solve phase is so dominant relative to the setup cost, the scaling limit for the end-to-end setup plus solve is dictated by the solve phase.
As such, the strong scaling limit occurs at $n/P\sim\num{27500}$, with the overall time being 1.28 seconds.
In this scenario, the resource usage is the same as the strong scaling limits identifier in \Cref{sec:liquid-phase-voltage,sec:solid-phase-voltage}.
For the liquid voltages in \Cref{sec:liquid-phase-voltage}, this translated to an end-to-end solve plus setup time of 0.4 seconds.
The end-to-end time for the solid voltages in \Cref{sec:solid-phase-voltage} is 0.48 seconds.
If the two voltages were uncoupled (i.e., the off-diagonal blocks in \Cref{eq:voltage-block} were zero), one would expect to solve the two voltages in 0.88 seconds, assuming no overlapping between the two solves.
While the relative residual tolerance considered here is two orders of magnitude larger than the tolerances in \Cref{sec:liquid-phase-voltage,sec:solid-phase-voltage}, this nevertheless demonstrates that the block Gauss-Seidel coupling in \Cref{eq:voltage-bgs} remains an effective preconditioner.

The weak scalability of the solver here remains comparable to the two subblock solvers in the two preceding subsections.
As shown in \Cref{tab:iteration-counts}, the coupled voltage system requires consistent iteration counts across problem scales, demonstrating the robustness of the block Gauss-Seidel preconditioner.
The setup phase exhibits an efficiency of 87\%.
The weak scalability of the solve phase is better at 93\%.
Finally, the end-to-end weak scalability is 91\%.
These numbers closely mirror that of the solid voltages in \Cref{sec:solid-phase-voltage}, whose cost exceeds the liquid-phase voltages from \Cref{sec:liquid-phase-voltage}.
One notable challenge with the solid voltages are the large jumps in the electrical conductivity in \Cref{eq:solid-phase-voltage}, which remain a focus for ongoing work on improving the strength-of-connection algorithms upon which SA-AMG techniques heavily rely.

\begin{table}[h]
    \caption{
      Iteration counts for subblock solvers and end-to-end electrochemical solve at the strong scaling limit (\(N/P=8\)).
      A relative residual tolerance of $10^{-8}$ is used for individual subblocks.
      The coupled voltage system, non-voltage system, and end-to-end solve use a $10^{-6}$ relative residual criterion.
      Results are shown for varying problem scales (\(N=1, 4, 16, 64\)), highlighting the performance of individual subblock preconditioners and the hierarchical block Gauss-Seidel preconditioner.
    }
    \centering
    \begin{tabular}{lcccc}
    \toprule
    Subblock & \(N=1\) & \(N=4\) & \(N=16\) & \(N=64\) \\ \hline
    \midrule
    Liquid-Phase Species (\Cref{sec:liquid-mass-fraction}) & 14 & 27 & 26 & 20 \\ \hline
    Liquid-Phase Pressure (\Cref{sec:liquid-phase-continuity}) & 24 & 29 & 31 & 29 \\ \hline
    Liquid-Phase Voltage (\Cref{sec:liquid-phase-voltage}) & 35 & 39 & 36 & 42 \\ \hline
    Solid-Phase Voltage (\Cref{sec:solid-phase-voltage}) & 32 & 32 & 33 & 32 \\ \hline
    \midrule
    Coupled Voltages (\Cref{sec:coupled-voltage}) & 32 & 28 & 29 & 28 \\ \hline
    Non-Voltage System (\Cref{eq:non-voltage}) & 27 & 37 & 42 & 38 \\ \hline
    \midrule
    End-to-End Solve (\Cref{sec:end-to-end-solver}) & 1 & 1 & 1 & 1 \\ \hline
    \bottomrule
    \end{tabular}
    \label{tab:iteration-counts}
\end{table}

\mysec{End-to-End Solver Scalability}\label{sec:end-to-end-solver}

\begin{figure}[h]
    \includegraphics[width=\textwidth]{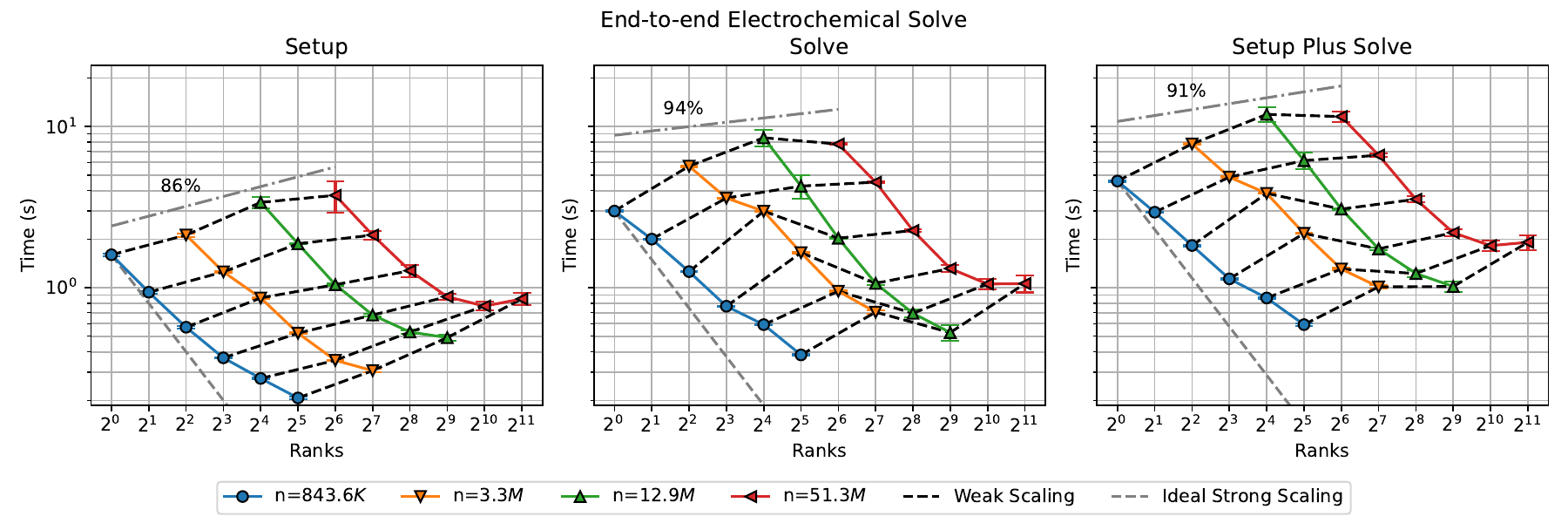}
    \caption{
        \label{fig:end-to-end-solve}
        Scalability of the end-to-end electrochemical solve.
        End-to-end setup and solve times reach 2.2 seconds for the largest problem scale ($n=51.3$ million degrees of freedom), with strong scalability observed up to $n/P\sim\num{100000}$.
        The weak scaling efficiency of 91\% highlights the scalability and efficiency of the hierarchical block preconditioning strategies for large-scale thermal battery simulations.
    }
\end{figure}

Armed with the subblock solvers identified in \Cref{sec:subblock-solvers},
we now look at the parallel scalability of the hierarchical block Gauss-Seidel preconditioner in \Cref{eq:hierarchical-bgs-exact-inverse}.
Flexible GMRES(30) is used to to apply the inverse of the non-voltage block, $\A_{n,n}^{-1}$, and inverse of the coupled voltage block, $\A_{v,v}^{-1}$.
Both systems are solved to a $10^{-6}$ relative residual tolerance.

The block Gauss-Seidel approach in \Cref{eq:non-voltage-bgs} is used to precondition the non-voltage block.
The liquid-phase species fraction subblock is approximated via a single sweep of the DD(0)-ILU(0) preconditioner from \Cref{sec:liquid-mass-fraction}.
Similarly, a single sweep of the the SA-AMG preconditioner from \Cref{sec:liquid-phase-continuity} approximates the solution to the liquid-phase pressure subblock.
For the coupled voltage solve, the block Gauss-Seidel preconditioner from \Cref{eq:voltage-bgs} in \Cref{sec:coupled-voltage}.
The inverse to each subblock is approximated via a single sweep of a SA-AMG preconditioner using the methods described in \Cref{sec:liquid-phase-voltage} and \Cref{sec:solid-phase-voltage} for the liquid and solid-phase voltages, respectively.
Results for the end-to-end electrochemical solve are shown in \Cref{fig:end-to-end-solve}.

Because the setup phase is executed sequentially across subblocks, the strong scaling of each subblock becomes a bottleneck for the monolithic preconditioner.
The three subblocks (\Cref{sec:liquid-phase-continuity,sec:liquid-phase-voltage,sec:solid-phase-voltage}) requiring construction of a SA-AMG preconditioner are the biggest source of bottlenecks.
Each of these reach their respective strong scale limits at the $P/N = 8$ level (the third data point from the left in the strong scaling results in \Cref{fig:end-to-end-solve}).
The highly scalable setup for the solid-phase conduction and liquid-phase species fraction block from \Cref{sec:liquid-mass-fraction} are sufficient to shift the strong scale limit to the right to $P/N = 16$.
This translates to a strong scale limit of $n/P\sim\num{100000}$ for the setup phase in the monolithic system,
which takes 0.88 seconds in the case of the largest scale considered.
However, the setup phase for each individual subblock preconditioner has no dependence on the setup of any other subblock.
A potential strategy not explored in this work is partitioning processors into disjoint subsets, with each subset dedicated to a specific subblock.
Then, the construction of the individual subblocks could be done in parallel.
As a result, the strong scale limit for the setup phase would be shifted to the right as this approach would exploit the embarassingly parallel nature of the setup phase.
Assuming little overhead associated with the partitioning and rebalancing approach, the strong scale limit for the setup phase then becomes limited by the least scalable subblock setup.
In this case, that would correpond to the strong scale limit of the solid-phase voltage block in \Cref{sec:solid-phase-voltage}, which is $n/P\sim\num{32000}$.
This represents the potential for a factor 3 improvement in the strong of the setup phase of the monolithic preconditioner, but would require knowledge about the strong scale limits of each individual subblock and partitioning the processors into disjoint subsets for each subblock.

The solve phase in \Cref{fig:end-to-end-solve} suffers from the same issue that the setup phase experiences.
In particular, the hierarchical block Gauss-Seidel preconditioner in \Cref{eq:hierarchical-bgs-exact-inverse} serializes the execution of the non-voltage and voltage solves due to the coupling through the $\A_{v,n}$ block.
Further, the solution to the $\A_{v,v}$ and $\A_{n,n}$ blocks is further serialized through the individual block Gauss-Seidel preconditioners used for each respective system.
The strong scaling limit here resembles that of the setup phase with $P/N=16$, albeit with better scaling at $\etastrong{}\sim0.74$ for the $N=64$ scale problem as compared to $\etastrong{}\sim0.53$ for the setup phase.
Nevertheless, this still translates to a relatively poor $n/P\sim\num{100000}$ at this limit.
As a result, the solution time is as large as 1.32 seconds.
A noticeable further improvement \emph{past} the strong scale limit may be made with $P/N=32$, leading to a solution time of 1.06 seconds.
The combined setup and solve scaling limit thus remains at $n/P\sim\num{100000}$ with an end-to-end time to solution of 2.2 seconds.
Finally, as shown in \Cref{tab:iteration-counts}, the end-to-end electrochemical solve converges in a single iteration by sub-iterating the voltage and non-voltage systems.
This demonstrates that the block-based preconditioner in \Cref{eq:hierarchical-bgs-exact-inverse} is robust.

While the setup phase may be further improved by considering exploiting the embarassingly parallel nature of the preconditioner construction of the various subblocks,
the same trick may not be applied here without exposing additional parallelism through altering the block Gauss-Seidel preconditioners in
the coupled voltage preconditioner in \Cref{eq:voltage-bgs}, coupled non-voltage preconditioner in \Cref{eq:non-voltage-bgs}, or even the hierarchical preconditioner in \Cref{eq:hierarchical-bgs-exact-inverse}.
One potential means of achieving this is foregoing the off-diagonal coupling term and instead employing a block Jacobi approach.
This approach would make the preconditioner evaluation embarrassingly parallel, allowing each subblock to be evaluated independently.
As before, this converts the inherently sequential action of the block Gauss-Seidel preconditioner to one that is instead bottlenecked by the slowest individual subblock solver.

The parallel block Jacobi approach described in the preceding paragraph is not without drawbacks, however.
Intuitively, this converts the block preconditioner from one that assumes there are non-trivial coupling terms between blocks to one in which the various physics blocks have no coupling.
As a consequence, the convergence rate of a block Jacobi method will suffer as a result.
One potential way of mitigating this is to employ a monolithic AMG approach, such as those considered by Ohm et al. \cite{ohm2021monolithic} or Fang et al. \cite{fang2019parallel}.
There, one could imagine a multigrd scheme wherein block Jacobi, which is embarassingly parallel with respect to the subblocks, is used as a \emph{smoother} for a multigrid approach.
In that case, the multiphysics couplings are still accounted for in the residual evaluation, as well as the coarser grid operators at each level.

Although several improvements could be made to the block Gauss-Seidel preconditioners proposed here, the strategies in \Cref{sec:block-preconditioning} remain effective in this context.
While a popular choice in many of \aria{}'s multiphysics applications, the one-level Schwarz schemes from \Cref{sec:domain-decomp} show poor convergence.
As a result solvers for the thermal batteries relied heavily on the sparse direct solvers provided through \klu{} prior to this work.
While robust, these solvers are nevertheless quite expensive.
While the symbolic LU factorization could be recycling between linear solves in an \aria{} simulation,
the numeric LU factorization would have to be repeated with every nonlinear iteration.
As a consequence, the bulk of the simulation time was occupied by the solve costs, which themselves were entirely dependent on the matrix factorization step.
This exhibits poor strong and weak scalability since \klu{} relies on gathering the matrix onto a single processor to perform the factorization.
For reference, the end-to-end direct solve cost for a factor 16 reduction in radial refinment to the baseline $N=1$ case took roughly 2 seconds on $P=16$ ranks.
Under the same resource utilization, the end-to-end setup plus solve cost for the $N=1$ baseline problem, which itself is a 16-fold radial refinement of the aforementioned case, takes only 0.86 seconds.
The approaches described in \Cref{sec:block-preconditioning} have enabled the solution to problems a factor 1024 times larger than the original, unrefined case.
Prior to this work, those problem scales would have been entirely intractible.

\mychp{Conclusion}
\label{chp:conclusion}
Effective and scalable preconditioning strategies are essential for achieving rapid time-to-solution in multiphysics simulations.
In the current work, we develop a highly scalable, block-based preconditioner for the simulation of thermal batteries.
We conduct weak scaling studies on all of the solver components required for the end-to-end solve, using several levels of mesh refinement up to $E=8.1M$ elements and $n=51.3M$ degrees of freedom.
Strong scaling is also observed, ranging from a single processor up to \num{2048} ranks.
We demonstrate that the solvers proposed here have good weak scalability with a combined 91\% parallel efficiency in the end-to-end electrochemical solve and strong scale to as many as \num{1024} ranks before providing no additional benefit to the time-to-solution.
As a result, we are able to achieve sub-second setup costs and nearly sub-second linear solves for end-to-end electrochemical solve.
While the results are currently limited to two dimensional axisymmetric calculation, they nevertheless demonstrate that the solver approaches discussed here scale well, and would also extend to full three-dimensional thermal battery models.

The strong scaling of the setup phase is limited by the sequential construction of individual subblock solvers, even though there are no dependencies between blocks.
We propose investigating means of parallelizing this setup across the subblocks in order to mitigate the strong scaling limits of the setup phase characterized by the preceding section.
We establish that the strong scaling of the block-based solvers is similarly fundamentally limited by the sequential nature of the hierarchical block Gauss-Seidel preconditioner.
Unfortunately, this indicates that the strong scaling limit for the solve phase is the sum of the scaling limits for each individual block.
In order to address these concerns, we propose further investigation into using block Jacobi preconditioners to repartition subblocks onto disjoint subsets of the processors.
This would allow several subblock evaluations to be done in parallel, extending the strong scaling limit further to the right.
One notable tradeoff with this proposed approach, however, is the impact on the preconditioner convergence by neglecting the coupling contribution in the strictly block upper triangular portion of the preconditioner.

Despite the areas for improvement noted above, this work represents a transformative advancement for thermal battery simulations.
As noted in the preceding section, \aria{} relied heavily on sparse direct solvers for thermal battery simulations, which exhibited poor scalability.
Through the methods developed in this work, we have extended the abilities of \aria{} to solve problems up to a factor \num{1024} times larger than were previously considered while providing considerable speedup to the linear solve phase.
As a result, TABS users are able to complete their calculations significantly faster, allowing for higher fidelity models to be considered.

\mychp*{Acknowledgments}
Sandia National Laboratories is a multi-mission laboratory managed and operated by National Technology \& Engineering Solutions of Sandia, LLC (NTESS), a wholly owned subsidiary of Honeywell International Inc., for the U.S. Department of Energy's National Nuclear Security Administration (DOE/NNSA) under contract DE-NA0003525.
This written work is authored by an employee of NTESS.
The employee, not NTESS, owns the right, title and interest in and to the written work and is responsible for its contents.
Any subjective views or opinions that might be expressed in the written work do not necessarily represent the views of the U.S. Government.
The publisher acknowledges that the U.S. Government retains a non-exclusive, paid-up, irrevocable, world-wide license to publish or reproduce the published form of this written work or allow others to do so, for U.S. Government purposes.
The DOE will provide public access to results of federally sponsored research in accordance with the DOE Public Access Plan.

\bibliographystyle{siamplain}
\bibliography{references}
\end{document}